\title{\large \bf A COMBINATORIAL METHOD
 FOR COMPUTING
STEENROD SQUARES}
\author{\large Roc\'{\i}o Gonz\'{a}lez-D\'{\i}az
and Pedro Real\\
Dpto. de Matem\'{a}tica Aplicada I \\
         Facultad de Inform\'{a}tica y
Estad\'{\i}stica \\
        Universidad de Sevilla\\
         rogodi@us.es,
         real@us.es}
  \date{ }

\documentstyle[emlines2,bezier,12pt]{article}
\newcommand{\Z}{{\bf Z}}

\newcommand{\ra}{\rightarrow}

\newcommand{\N}{\bf N}

\newcommand{\scst}{\scriptscriptstyle}
\newcommand{\ot}{\otimes}
\newtheorem{prp}{Proposition}[section]
\newtheorem{lmm}[prp]{Lemma}
\newtheorem{thr}[prp]{Theorem}

\newtheorem{crl}[prp]{Corollary}

\begin{document}

\bibliographystyle{alpha}

\sloppy

\maketitle

\vspace{0.5cm}

\begin{center}
{\bf Abstract}
\end{center}
We present here a combinatorial method for computing cup-$i$
products and
Steenrod squares of a simplicial set $X$. This method
is essentially  based on the determination of  explicit formulae for the
component morphisms of a higher diagonal approximation
(i.e., a family of morphisms measuring the lack of commutativity of
the cup product on the cochain level)
in terms of  face operators of $X$.
A generalization of this method to Steenrod reduced
powers is sketched.

\vfill

\noindent {\scriptsize Both authors are partially supported by the
PAICYT research project
 FQM-0143 from Junta de Andaluc\'{\i}a and the DGESIC research project
PB97-1025-C02-02 from Education and Science Ministry (Spain).}

\newpage

\section{Introduction}

Cohomology operations are algebraic operations on the cohomology
groups of spaces commuting with the homomorphisms induced by
continuous mappings. This machinery is useful when the graded
vector space structure and the cup product on cohomology fail
to distinguish two spaces by their cohomology. Steenrod squares
\cite{steenrodepstein62}
constitute an extremely important class of cohomology operations
not only in Algebraic Topology but also in the area of simplicial
methods in Homological Algebra (cohomology of groups, Hochschild
cohomology of algebras, ...).

 We here study in detail the underlying
combinatorial structures of the definition of these
cohomology operations. More concretely, we will move in the framework of
Simplicial Topology \cite{may67} in which the basic objects are
simplicial sets, that is, graded sets endowed with face
($\partial_i$) and degeneracy ($s_i$) operators, satisfying
several commutativity relations. Roughly speaking, a simplicial
set can be considered as an algebraic generalization of the
structure of a triangulated polyhedron although the former features a
more rigid combinatorial structure than the latter. In this
context, we develop a suitable setting in which all Steenrod
cohomology operations can be studied simultaneously.

It is well-known that there are  several methods
for constructing Steenrod squares. One of them consists
of making these operations using the cohomology
of Eilenberg-Mac Lane spaces (see, for instance, [\cite{may67}; p. 107]).
Another method is determined by the construction of a family
of morphisms $\{D_i\}$ measuring   the lack  of
commutativity of the cup product on  the cochain level
\cite{steenrodepstein62}. This sequence of morphisms is
called {\em higher diagonal
approximation}. Its existence is always guaranteed by
the acyclic models method
\cite{eilenbergmaclane53}. In this way, it is possible to derive a
recursive procedure to obtain the explicit formula for any
$D_i$ (see \cite{DOLD}, [\cite{may70}; Sect. 7]).

In this paper, we present an alternative method for
obtaining the explicit formula for a higher diagonal approximation.
 In \cite{real}, the formula for
$D_i$ is established in terms of the
component morphisms of a given
Eilenberg-Zilber contraction (a special homotopy equivalence)
 from $C_{*}^{\scst N}(X\times X)$
onto $C_{*}^{\scst N}(X)\ot
C_{*}^{\scst N}(X)$, where $C_{*}^{\scst N}(X)$ denotes the
normalized chain complex  of a given simplicial set $X$.

 Now, we have  the following problems.  On one hand, the  component
morphisms of the contraction above are defined in terms of
face  and degeneracy operators
 of the simplicial set $X$. On the other hand,
the formula for a morphism $D_i$
always involves the use of the homotopy
operator of the Eilenberg-Zilber contraction
and the explicit formula for this last morphism is
determined by shuffles (a special type of permutation) of degeneracy
operators. In consequence, if we try to express in this way  the
morphisms $D_i$ in terms of face and degeneracy operators of $X$,
the number of summands  appearing in
the formula for a morphism $D_i$ evaluated over an
element of degree $n$ is, in general, at least $2^n$.
Therefore, an algorithm that would be designed starting from these
formulae  would be too
slow for practical implementation.

Because of this, the idea of
simplifying these formulae arises in a natural way.
This simplification or normalization  is based on the fact that
any composition of
face and degeneracy operators of the simplicial set $X$ can be put
in a ``canonical'' way. That is, a composition
of this type can be expressed   in the unique form:
$$s_{j_t}\cdots s_{j_1}\partial_{i_1}\cdots
\partial_{i_s},$$
where  $j_t>\cdots >j_1\geq 0$
 and $i_s>\cdots >i_1\geq 0$.

Moreover, taking into
account that the image of a morphism $D_i$ lies
in \mbox{$C_{*}^{\scst N}(X)\ot C_{*}^{\scst N}(X)$},
  those summands of the simplified formula for
$D_i$  with a
factor having a degeneracy
operator in its expression, can be eliminated. The
reason is that this factor applied to an element
of the simplicial set $X$ is zero in the normalized
chain complex associated to $X$. In this way, we here
obtain more simple  formula for the morphism
$D_i$.

In this way, the classical definition of Steenrod squares:

\begin{eqnarray}\label{SQ}
Sq^i(c)(x)&=&\left\{\begin{array}{ll}
\mu(<c\ot c, \; D_{j-i}(x)>),\; & i\leq j \\
0, & i>j \end{array} \right.  \end{eqnarray}
where $c\in \mbox{Hom}(C_j^{\scst N}(X), \Z_2)$,
$x\in C_{i+j}^{\scst N}(X)$ and $\mu$ is the homomorphism induced
by the multiplication in $\Z_2$,
is complemented by a manageable combinatorial
formulation of  the higher
diagonal approximation. As we will remark later, this
description can be considered as a direct translation of the
most ancient  definition of Steenrod squares (see \cite{steenrod47})
to the  general setting of the Simplicial Topology.
We give a more detailed explanation in the third section.

 Still, we think that this combinatorial machinery could be
   substantially improved in the future by exploiting the well-known
properties
   of Steenrod squares \cite{adem52,steenrodepstein62} and advanced
   techniques  for calculating cocycles (see, for example,
\cite{egl97,lambe97}). In this way, a ``reasonably
efficient'' algorithm
computing the cohomology algebra of several important simplicial
sets could be derived.

Finally, we start an analogous  study  to that given  in
\cite{real} for Steenrod reduced  powers. More precisely,
we provide  a simplicial  description of these operations
in terms of the component morphisms of
an Eilenberg-Zilber contraction from
$C_{*}^{\scst N}(X\times \stackrel{\mbox{\scriptsize $p$
times}}\cdots \times X)$ to $ C_{\star}^{\scst N}(X)\ot
\stackrel{\mbox{\scriptsize $p$
times}}\cdots \ot C_{*}^{\scst N}(X)$,
and  a $0$-sequence $\{\gamma_i\}_{i \geq 0}$ in the symmetric
group $G_p$.

Here is a summary of the present paper. Section 2 is dedicated to notation,
terminology and a presentation of the problem.
In Section 3, we show an explicit
combinatorial definition of cup-$i$ products and, consequently, of
Steenrod squares. In Section 4, we
study Steenrod reduced powers from this point of view.
Finally, a proof of the main theorem enunciated in Section 3 is
given in Section 5.

We are grateful to Prof. Tom\'{a}s Recio for his helpful suggestions for
making the exposition more readable and, hopefully, more informative.
We also wish to acknowledge our debt to the referees for their many
valuable indications.

\section{The higher homotopy commutativity of the
Alexander-Whitney operator}

  The aim of this section is to give some preliminaries and a
  brief account of the work done in \cite{real} in order to
  facilitate the understanding of the rest of the paper. Most of the
material given in this section can be found in
  \cite{may67}, \cite{maclane75}, \cite{spanier} and \cite{weibel}.

 A {\em simplicial set} $X$ is a sequence of sets $X_0, X_1,
 \ldots$, together with {\em face operators} $\partial_i:
 X_n \rightarrow X_{n-1}$ and {\em degeneracy operators}
 $s_i: X_n \rightarrow X_{n+1}$ ($i=0,1,\ldots, n)$, which
 satisfy the following ``simplicial'' identities:

\begin{eqnarray*}
\mbox{\bf (s1)} & & \partial_i \partial_j = \partial_{j-1} \partial_i, \quad
\mbox{ if } i<j;\\
\mbox{\bf (s2)} & & s_i s_j= s_{j+1} s_i,\quad \mbox{ if } i \leq j;\\
\mbox{\bf (s3)} & & \partial_i s_j = s_{j-1} \partial_i,\quad \mbox{ if }
i<j,\\
\mbox{\bf (s4)}& &\partial_i s_j = s_j \partial_{i-1},\quad
\mbox{ if } i>j+1,\\
\mbox{\bf (s5)}& & \partial_j s_j = 1_{\scst X} =
\partial_{j+1} s_j.\end{eqnarray*}

  The elements of $X_n$ are called $n$-{\em simplices}. A simplex $x$ is
  {\em degenerated} if $x=s_i y$ for some simplex $y$ and
  degeneracy operator $s_i$; otherwise, $x$ is non degenerated.

  The following elementary
  lemma will be essential in the proof of the main
  \mbox{theorem} of this paper.

\begin{lmm}\cite{may67} \label{can}
Any composition $\mu: X_m \rightarrow X_n$ of face and degeneracy
operators of a simplicial set $X$ can be put in a unique
``canonical'' form:

\begin{eqnarray*}
s_{j_t}\cdots s_{j_1}\partial_{i_1}\cdots
\partial_{i_s},
\end{eqnarray*}
where  $n>j_t>\cdots >j_1\geq 0$,
 $m\geq i_s>\cdots >i_1\geq 0$ and $n-t+s =m$.
\end{lmm}

Let $R$ be a ring which is commutative with unit. A {\em chain}
(resp. {\em cochain}) {\em complex}
$M=\{M_n, d_n\}$ (resp. $C=\{C^n, \delta^n\}$) is a graded (over the
integers) \mbox{$R$-module} together with a $R$-module map
$d$ of degree $-1$ (resp. of degree $+1$), called the {\em differential},
such that
$d_n \, d_{n+1} = 0$ (resp. $\delta^{n+1} \, \delta^n =0$). An
element $c$ of $C^n$ is called $n$-cochain.
The {\em homology} $H_{*}(M)$ (resp.  the {\em cohomology}
$H^{*} (C)$ is the family of modules

$$H_n (M) \; = \; \mbox{\rm Ker } d_n / \mbox{\rm Im } d_{n+1}
\;\;\;\;\;\;\; \mbox{(resp. }\; H^n (C) =
\mbox{\rm Ker } \delta^{n} / \mbox{\rm Im } \delta^{n-1}).$$

 If $M$ is a chain complex over $R$ and $G$ is an $R$-module,
 there is a cochain complex $\mbox{Hom($M,G$)}=
 \{\mbox{Hom($M_n, G$)}, \delta^n\}$, where, if $c\in
 \mbox{Hom($M_n, G$)}$, then \mbox{$\delta^n c \in \mbox{Hom($M_{n+1},
 G$)}$} is defined by

$$(\delta^n c)(x) = c (d_{n+1} (x)),\;\;\;\;\;\;\; x\in M_{n+1}.$$

 We also write $<c, x>$ instead of $c(x)$ and set $<c,x>=0$ if
 the degree of the cochain $c$
 is not equal to the degree of the element $x$. In this notation,

$$ < \delta^n c, x> \; = \; <c, d_{n+1} (x)>. $$

  Now, given a simplicial set $X$, let $C_*(X)$
  denotes the chain complex
$\{ C_{n} (X), d_n\}$, in which $C_n(X)$ is the
  free
  $R$-module generated by $X_n$ and \mbox{$d_n : C_n(X) \rightarrow
  C_{n-1}(X)$} is defined by $\displaystyle
  d_n = \sum_{i=0}^{n} (-1)^{i} \partial_i$. Let us denote by
  $s(C_* (X))$ the graded \mbox{$R$-module} gene\-rated by all the
  degenerated simplices.  In $C_*(X)$, we have that
  \mbox{$d_n(s(C_{n-1}(X))) \subset s(C_{n-2}(X))$} and, then
  $C_{*}^{\scst N} (X)= \{ C_n (X)/s(C_{n-1} (X)), d_n\}$ is a
  chain complex called {\em the normalized chain complex
  associated
  to $X$}. Given a \mbox{$R$-module} $G$, let us define by $C^*(X;G)$ the
  cochain complex associated to $C_*^{\scst N}(X)$. In this way, we
  define the homology and cohomology of $X$
  with coefficients in a $R$-module $G$ by
  $H_*(X; G) = H_*(C_*^{\scst N}(X) \otimes G)$ and
  $H^{*}(X; G) = H^{*}(C^*(X;G))$,
  res\-pectively.

Eilenberg and Mac Lane defined in \cite{eilenbergmac53}
a {\em contraction} of chain complexes
from $N$ onto $M$, as a triple $(f,g,\phi)$ in
which $f:N\ra M$ (projection) and $g:M\ra N$ (inclusion)
are chain maps, and
$\phi:N\ra N$ (homotopy operator) is
 a map of $R$-module raising degree by $1$. Moreover, it
is required that
$$\mbox{\bf (c1)} \quad fg= 1_{\scst M},
\qquad\mbox{\bf (c2)}\quad \phi d + d\phi + gf = 1_{\scst N},$$
$$ \mbox{\bf (c3)} \quad \phi g = 0,
\qquad\mbox{\bf (c4)} \quad  f\phi = 0,
\qquad\mbox{\bf (c5)} \quad  \phi\phi = 0.$$

   Hence, this definition implies that the ``big'' complex $N$ is
   homology equivalent to the ``small'' complex $M$ in a strong
   way. In fact, a contraction is a special homotopy equivalence
   between chain complexes.

If we have two contractions  $(f_i, g_i,
\phi_i)$ from $N_i$ to $M_i$, with $i=1,2$ then,  the following
contractions can be constructed (see
\cite{eilenbergmac53}):
\begin{itemize}
\item The tensor product contraction
$(f_1\ot f_2, g_1\ot g_2, \phi_1\ot g_2 f_2+1_{\scst N_1}\ot
\phi_2)$ from $ N_1\ot N_2$ to $M_1\ot M_2$.

\item If $N_2=M_1$, the  composition contraction
$(f_2 f_1, g_1 g_2, \phi_1+g_1\phi_2 f_1)$ \mbox{from $ N_1$} to $ M_2$.
\end{itemize}

    If $p$ and $q$ are non-negative integers, a {\em
    $(p,q)$-shuffle} $(\alpha, \beta)$ is a partition of the set
    $\{0,1, \ldots, p+q-1\}$ of integers into two disjoint
    subsets, $\alpha_1 < \cdots < \alpha_p$ and
    $\beta_1 < \cdots < \beta_q$, of $p $ and $q$ integers,
    respectively. The signature  of the
    shuffle $(\alpha, \beta)$ is defined by
    $sig(\alpha, \beta)= \displaystyle\sum_{i=1}^{p} \alpha_i - (i-1).$

 After these preliminaries, we are
  able to describe a very important homotopy equivalence
   in Algebraic Topology. This contraction tells us that the
   associated chain complex $C_*^{\scst N}(X\times Y)$ reduces to
   the
   tensor product of chain complexes $C_*^{\scst N} (X)$ and
   $C_*^{\scst N}(Y)$.

  An Eilenberg-Zilber contraction \cite{eilenbergzilber59} from
$C_{*}^{\scst N}(X\times Y)$ to $C_{*}^{\scst N}(X)\otimes
C_{*}^{\scst N}(Y)$,
where $X$ and $Y$ are given simplicial sets, is defined by the
triple $(AW,\,EML,\,SHI)$ where:

\vspace{1cm}

\begin{itemize}
\item The Alexander-Whitney operator
$AW:\,C_{*}^{\scst N}(X \times Y) \longrightarrow
C_{*}^{\scst N}(X)\ot C_{*}^{\scst N}(Y)$ is defined by:

$$AW(a_m \times b_m)= \sum _{i=0}^m \partial _{i+1} \cdots \partial_m
a_m \ot \partial _0 \cdots \partial _{i-1} b_m.$$

  If $X=Y$, $AW$ can be considered as a ``simplicial approximation'' to the
  dia\-gonal and this operator provides a method for constructing the
  cup product in
  cohomology.  If we interchange the factors $a_m$ and
  $b_m$ in the formula, we obtain a different approximation.
  Comparison of these two different approximations to the diagonal
  leads to the Steenrod squares.

\item The Eilenberg-Mac Lane operator
$EML:\, C_{*}^{\scst N}(X)\ot C_{*}^{\scst N}(Y)\longrightarrow
C_{*}^{\scst N}(X \times Y)$ is defined by:

$$EML(a_p \ot b_q) = \sum _{( \alpha , \beta ) \in \{ (p,q)-
\mbox{\scriptsize shuffles} \}} (-1)^{sig( \alpha , \beta )}s_{\beta _q}
\cdots s_{\beta _1}a_p \times s_{\alpha _p} \cdots
s_{\alpha_1}b_q.$$

    This operator can be seen as a process of ``triangulation'' in
    the cartesian product $X\times Y$.

\item And the Shih operator
$SHI:\,C_{*}^{\scst N}(X \times Y)\longrightarrow
C_{*+1}^{\scst N}(X \times Y)$ is defined by:

$$\begin{array}{lll}
SHI(a_0\times b_0)&=0;&\\\\
SHI(a_m \times b_m) &=
\displaystyle\sum (-1)^{\bar{m}+sig( \alpha , \beta )+1} &s_{\beta
_q+\bar{m}} \cdots s_{\beta _1+\bar{m}}s_{\bar{m}-1} \partial _{m-q+1}
\cdots \partial _m a_m\\
&& \times s_{\alpha _{p+1}+\bar{m}} \cdots
s_{\alpha_1+\bar{m}} \partial _{\bar{m}} \cdots \partial _{m-q-1}b_m;
\end{array}$$
\noindent where $\bar{m}=m-p-q$, $sig( \alpha , \beta )=
\displaystyle\sum _{i=1}^{p+1} \alpha _i
-(i-1)$, and the last sum is taken over all the indices
$0 \leq q \leq m-1$, $0\leq p \leq m-q-1$ and
\mbox{$( \alpha , \beta ) \in \{ (p+1,q)$-shuffles$\}$.}
\end{itemize}

  A recursive formula for the $SHI$ operator has already been given by
  Eilenberg and Mac Lane in \cite{EMII}. The explicit formula
  given here was stated by Rubio in
 \cite{rubio91} and proved by Morace in the appendix of \cite{real96}.
 In contrast to the deep studies found in the literature  on
 the $AW$ and $EML$ operators, it turns to be  quite surprising
 the lack of interest
 shown up to now in the study of the homotopy operator involved in an
 Eilenberg-Zilber
contraction, not only from the point of view of getting its explicit
formula but also of obtaining algebraic preservation results of
this
operator with regard to the underlying coalgebra
structure on
$C_*^{\scst N}(X\times X)$.

Given a simplicial set $X$ and a positive integer
$p$, we can form  a contraction
$(f,g,\phi)$ from $ C_{*}^{\scst N}(X\times
\stackrel{\mbox{\scriptsize $p$ times}}\cdots \times X)
$ to  $ C_{*}^{\scst N}(X)\ot \stackrel{\mbox{\scriptsize $p$
times}}\cdots \ot C_{*}^{\scst N}(X)$,
appro\-priately composing  Eilenberg-Zilber contractions. If
$p=2$, then $(f,g,\phi)$ is the contraction $(AW, EML,
SHI)$.  The contraction from
$C_{*}^{\scst N}(X\times X\times X)$ onto \mbox{$C_{*}^{\scst N}(X)\ot
C_{*}^{\scst N}(X)\ot C_{*}^{\scst N}(X)$} is defined by
the composition of the
Eilenberg-Zilber contraction from $C_{*}^{\scst N}(X\times
X\times X)$ to $C_{*}^{\scst N}(X)\ot C_{*}^{\scst N}
(X\times X)$ and the tensor product contraction
of the identity morphism $1_{C_{*}^{\scst N}(X)}$
 and the Eilenberg-Zilber contraction from
$C_{*}^{\scst N}(X\times X)$ onto $C_{*}^{\scst N}(X)\ot
C_{*}^{\scst N}(X)$. And so on.

  From now on, the contractions obtained in this way
  will also be called
Eilenberg-Zilber contractions.

Let $p$ be a positive integer. Let us define
several chain maps (or morphisms) we
 use in this paper. We will omit in the notation of these morphism
 its dependency on $p$. The {\em diagonal map}
$$\Delta: C_{*}^{\scst N}(X)\ra C_{*}^{\scst N}(X\times
\stackrel{\mbox{\scriptsize $p$ times}}\cdots \times X)$$
is defined by $\Delta(x)=(x,\stackrel{\mbox{\scriptsize
$p$ times}}\dots,x)$. The following  automorphisms are also defined:
$$t: C_{*}^{\scst N}(X\times
\stackrel{\mbox{\scriptsize $p$ times}}\cdots \times X)\ra
C_{*}^{\scst N}(X\times
\stackrel{\mbox{\scriptsize $p$ times}}\cdots \times
X),$$ such that $t(x_1, x_2, \dots, x_p)=(x_2, \dots,
x_p, x_1)$ and
$$T: C_{*}^{\scst N}(X)\otimes
\stackrel{\mbox{\scriptsize $p$ times}}\cdots \otimes
C_{*}^{\scst N}(X)\ra C_{*}^{\scst N}(X)\otimes
\stackrel{\mbox{\scriptsize $p$ times}}\cdots \otimes
C_{*}^{\scst N}(X)$$
defined by $T(x_1\ot x_2\ot \cdots
\ot x_p)=(-1)^{|x_1|(|x_2|+ \cdots +|x_p|)}
x_2\ot \cdots \ot x_p\ot x_1$.

\vspace{0.3cm}

  We now outline the problem concerning Steenrod squares in which
  we are inte\-rested. It is well-known that the $AW$ operator
is not commutative, that is, assuming that $X=Y$,
$AW t \neq AW$. On the other hand, this operator
determines the cup product in cohomology. If $G$ is a ring,
given two cochains
$c\in C^{i}(X;G)$ and $c'\in C^{j}(X;G)$, and
$x\in C_{i+j}^{\scst N} (X)$, the cup-product of $c$ and $c'$ is
defined by:

   $$\begin{array}{rcl}
   c \smile c' (x) & = & \mu ( <c \otimes c', AW\; \Delta (x)>) \\
                  &   & \\
                  &=  &
                  \mu (<c, \partial_{i+1} \cdots \partial_{i+j}x>
                  \ot <c', \partial_{0} \cdots
                  \partial_{i-1}x>),
                  \end{array}$$
where $\mu$ is the homomorphism induced by the multiplication on $G$.
\mbox{Steenrod} in \cite{steenrod52} determined that  there
exists an infinite sequence
of morphisms $\{D_i\}$, called {\em higher dia\-gonal
approximation},
which ``measures'' this lack of
commutativity. More precisely, there is a sequence of graded
homomorphisms \mbox{$D_i: C_*^{\scst N}(X) \ra C_*^{\scst N}(X) \ot
C_*^{\scst N}(X)$} of degree $i$ such that:

$$ \begin{array}{c}
D_0 \; = \; AW \; \Delta \\
d_{\scst\ot} \; D_{i+1} + (-1)^{i} D_{i+1} \; d \; = \;
T \; D_i + (-1)^{i+1} D_i,
\end{array}
$$
where $d$ and $d_{\scst \ot}$ are the differentials
of $C_*^{\scst N}(X)$ and $C_*^{\scst N}(X) \ot C_*^{\scst
N}(X)$, respectively.

  Moreover, the morphism
$D_i$ can be expressed in the form
$D_i=h_i\; \Delta$, where $\mbox{$h_i:
C_*^{\scst N}(X\times X) \rightarrow
C_*^{\scst N}(X)\otimes C_*^{\scst N}(X)$}$ is a
homomorphism of degree $i$. In the
lite\-rature, the
existence of the tower of iterated Alexander-Whitney operators
$\{h_i\}$ is guaranteed by
the acyclic models method
(see \cite{eilenbergmaclane53}). This technique can be conside\-red
as
a constructive method in the simplicial category and
recursive formulae for
$\{D_i\}$ can be established (see \cite{may70}).

In \cite{real}, a different approach is presented. The
strong homotopy commutativity of $AW$ is determined by making use of the
explicit Eilenberg-Zilber contraction $(AW, EML, SHI)$. Hence, the
formula for a morphism $D_i$ is given in
terms of the morphisms $AW$, $SHI$, the diagonal $\Delta$ and
the automorphism $t$. More precisely, the formula for
$h_i$ is $AW\;(t\; SHI)^i$, for all $i\in N$.

Now, the definition of the cohomology operation
$Sq^i:\,H^j(X; \Z_2)\ra H^{j+i}(X; \Z_2)$ (see (\ref{SQ})) takes the
form:

\begin{eqnarray}\label{1}
Sq^i(c)(x)=\left\{\begin{array}{ll}
\mu(<c\ot c, AW\; (t\; SHI)^{j-i}(x,x)>),\; & i\leq j \\
 0, & i>j \end{array} \right.
\end{eqnarray}
\noindent where $c\in C^j(X; \Z_2)$ and
$x\in C_{i+j}^{\scst N}(X)$.

 It is obvious that the formulae for the morphism $h_i$ can be
 given in terms of face and
degeneracy  operators of $X$. Our objective in the next section
is to show how to ``simplify''
these formulae and to obtain an explicit definition of Steenrod
squares only in terms of  face operators of $X$.

\section{An explicit combinatorial description of the
cup-$i$ products}

It is clear that the image of $h_i$ lies in
$C_*^{\scst N}(X)\ot C_*^{\scst N}(X)$.
Therefore, if we express the factors of the summands of
the formula for $D_i$ in a canonical way
(see Lemma \ref{can}), those summands of the
simplified formula for $D_i$
having a factor with a degeneracy operator in
its expression must be eliminated.

Working in this way, we can state
the following theorem. Section 5 is entirely devoted to the proof
of this result.

\begin{thr}\label{2}
Let $R$ be the ground ring and let $X$ be a
simplicial set. Let us consider the Eilenberg-Zilber
contraction $(AW, EML, SHI)$ from
$C_{*}^{\scst N}(X\times X)$ onto
$C_{*}^{\scst N}(X)\ot C_{*}^{\scst N}(X)$. Then,   the morphism
$h_n=AW(t\, SHI)^n:\,C_m^{\scst N}(X\times X)\ra
(C_{*}^{\scst N}(X)\ot C_{*}^{\scst N}(X))_{m+n}$
can be expressed in the form:

\begin{itemize}
\item if $n$ is even, then:
\begin{eqnarray*}
AW(t\, SHI)^n & = & \sum_{i_n=n}^{m}\quad
\sum_{i_{n-1}=n-1}^{i_n-1}\, \cdots \,
\sum_{i_0=0}^{i_1-1} \, (-1)^{A(n)+B(n,m,\bar{\i})+C(n,\bar{\i})+D(n,m,
\bar{\i})}\\ \\
& & \partial_{i_0+1}\cdots \partial_{i_1-1}\partial_{i_2+1}\cdots
\cdot \partial_{i_{n-1}-1} \partial_{i_{n}+1} \cdots
\partial_{m} \\
& & \ot \partial_0 \cdots \partial_{i_0-1} \partial_{i_1+1} \cdots
\cdot \partial_{i_{n-2}-1} \partial_{i_{n-1}+1} \cdots
\partial_{i_n-1}, \end{eqnarray*}

\item if $n$ is odd, then:
\begin{eqnarray*}
AW(t\,SHI)^n &=& \sum_{i_n=n}^{m}\quad
\sum_{i_{n-1}=n-1}^{i_n-1}\, \cdots \,
\sum_{i_0=0}^{i_1-1} \ (-1)^{A(n)+B(n,m,\bar{\i})+C(n,\bar{\i})+
D(n,m,\bar{\i})}\\\\
& & \partial_{i_0+1}\cdots \partial_{i_1-1}\partial_{i_2+1}\cdots
\cdot \partial_{i_{n-2}-1} \partial_{i_{n-1}+1} \cdots
\partial_{i_{n}-1}  \\
& & \ot \partial_0 \cdots \partial_{i_0-1} \partial_{i_1+1} \cdots
\cdot \partial_{i_{n-1}-1} \partial_{i_n+1} \cdots
\partial_m,
\end{eqnarray*}
\end{itemize}
\noindent where
\begin{eqnarray*}
A(n)&=&\left\{\begin{array}{l}
 1, \quad \mbox{ if } n\equiv 3,4,5,6 \mbox{ mod } 8\\
 0, \quad \mbox{ otherwise } \end{array}\right.\\\\
B(n,m, \bar{\i})&=&\left\{\begin{array}{ll}
\displaystyle\sum_{j=0}^{\lfloor\frac{n}{2}\rfloor}i_{2j}, & \mbox{
if } n\equiv 1,2\mbox{ mod } 4\\\\
\displaystyle\sum_{j=0}^{\lfloor\frac{n-1}{2}\rfloor}i_{2j+1}+n m, & \mbox{ if }
n\equiv 0,3 \mbox{ mod } 4
\end{array}\right.\\\\
C(n, \bar{\i})&=&\sum_{j=1}^{\lfloor\frac{n}{2}\rfloor}
(i_{2j}+i_{2j-1})(i_{2j-1}+\cdots+i_0)
\end{eqnarray*}
and
\begin{eqnarray*}
D(n,m, \bar{\i})&=&\left\{\begin{array}{ll}
(m+i_n)(i_n+\cdots +i_0),\quad &
\mbox{ if } n \mbox{ is odd }\\
 0,\quad & \mbox{ if } n \mbox{ is even}
\end{array}\right.
\end{eqnarray*}
where $\bar{\i}= (i_0, i_1, \ldots, i_{n})$.
\end{thr}

  A first consequence of the theorem above is that
taking into
account that the cup-$i$ product (see, for example,
\cite{dieudonne89})
of a $p$-cochain $c\in C^{p}(X;G)$
and a $q$-cochain $c'\in C^{q}(X;G)$ is defined by:

$$c\smile_i c'(x) \; = \; \mu (<c \ot c', D_i (x)>),
$$ where $i\in \N$ and $x \in C_{p+q-i}^{\scst N} (X)$, one immediately
gets a simplicial description of these operations. In an
analogous way, a combinatorial
definition of Steenrod squares (\ref{1}) is given by the following
corollary:

\begin{crl} \label{3}
 Let $\Z_2$ be the ground ring and let
$X$ be a simplicial set.
If \mbox{$c\in C^j(X; \Z_2)$} and
$x\in C_{i+j}^{\scst N}(X)$, then
$Sq^i:\,H^j (X; \Z_2)\ra
H^{j+i}(X; \Z_2)$ is defined by:

\begin{itemize}
\item If $i\leq j$ and $i+j$ is even, then:

\begin{eqnarray*}
Sq^i(c)(x) &=& \sum_{i_{n}=S(n)}^{m}\quad
\sum_{i_{n-1}=S(n-1)}^{i_{n}-1}\,\cdots \,
\sum_{i_1=S(1)}^{i_2-1}
\\ \\
& & \mu(<c,\partial_{i_0+1}\cdots \partial_{i_1-1}\partial_{i_2+1}\cdots
\cdot \partial_{i_{n-1}-1} \partial_{i_{n}+1} \cdots \partial_{m} x>
\\
& & \ot <c, \partial_0 \cdots \partial_{i_0-1} \partial_{i_1+1} \cdots
\cdot \partial_{i_{n-2}-1} \partial_{i_{n-1}+1} \cdots \partial_{i_n-1}
x>).
\end{eqnarray*}

\item If $i\leq j$ and $i+j$ is odd, then:

\begin{eqnarray*}
Sq^i(c)(x) &=& \sum_{i_{n}=S(n)}^{m}\quad
\sum_{i_{n-1}=S(n-1)}^{i_{n}-1}\,\cdots \,
\sum_{i_1=S(1)}^{i_2-1}\\ \\
& & \mu(<c,\partial_{i_0+1}\cdots \partial_{i_1-1}\partial_{i_2+1}\cdots
\cdot \partial_{i_{n-2}-1} \partial_{i_{n-1}+1} \cdots
\partial_{i_{n}-1} x>  \\
& & \ot <c, \partial_0 \cdots \partial_{i_0-1} \partial_{i_1+1} \cdots
\cdot \partial_{i_{n-1}-1} \partial_{i_n+1} \cdots
\partial_m x>).
\end{eqnarray*}

\item If $i>j$, then
$Sq^i(c)(x)=0.$
\end{itemize}

In these formulae, $n=j-i$, $m=i+j$,
$$S(k)=i_{k+1}-i_{k+2}+\cdots
+(-1)^{k+n-1}i_n+(-1)^{k+n} \lfloor
\frac{m+1}{2}\rfloor +\lfloor\frac{k}{2}\rfloor,$$ for
all $0\leq k\leq n$ and $i_0=S(0)$.
\end{crl}

\noindent{\bf Proof}

Let us start with $c\in C^j(X; \Z_2)$.
If $i>j$ then
$Sq^i(c)=0$. So, let us suppose that
$i\leq j$.
Then,
$$Sq^i(c)(x)=\mu(<c\ot c,\;AW(t\, SHI)^n(x,x)>)$$
where $x\in C_m(X)$.
It is not hard to notice that
we only have to consider the summands of
the explicit formula for $AW\,(t\, SHI)^n$
(see Theorem \ref{2}),  having the same
number of face operators in both factors.

\begin{itemize}
\item If $n=0$, then $m-i_0=i_0$, so
$i_0=\frac{m}{2}$. The
formula only has  one summand.

\item If $n=1$, then
$i_1-1-i_0=i_0+m-i_1$. So,
$i_0=i_1-\frac{m+1}{2}$
and
$i_1\geq \frac{m+1}{2}$.

\item In general, if $n$ is even (when $n$ is odd the
proof is analogous):
\begin{eqnarray*}
&&m-i_n+\cdots+i_{2k+1}-1-i_{2k}+\cdots+i_1-1-i_0  \\
&&=i_n-1-i_{n-1}+\cdots+i_{2k}-1-i_{2k-1}+\cdots+
i_2-1-i_1+i_0,
\end{eqnarray*}
hence, we have
\begin{eqnarray} \label{5}
i_0=i_1-i_2+i_3-\cdots -i_n+\frac{m}{2}.
\end{eqnarray}

Taking into account in (\ref{5}) that $i_0\geq 0$, we get
\begin{eqnarray*}
i_1\geq i_2-i_3+\cdots +i_n-\frac{m}{2}.
\end{eqnarray*}

Using $i_0\leq i_1-1$ in (\ref{5}), we have
\begin{eqnarray*}
i_2\geq i_3-i_4+\cdots +i_{n-1}-i_n+\frac{m}{2}+1.
\end{eqnarray*}

In general, let us suppose that
\begin{eqnarray*}
i_k\geq i_{k+1}-i_{k+2}+\cdots +
(-1)^{k+n-1} i_n + (-1)^{k+n}
\frac{m}{2}+\lfloor\frac{k}{2}\rfloor
\end{eqnarray*}
for all $1\leq k\leq \ell$ and let us prove that this
expression is true in $\ell+1$.
In the case $k=\ell-1$, since $i_{\ell}-1\geq i_{\ell-1}$, we have
\begin{eqnarray*}
i_{\ell}-1\geq i_{\ell}-i_{\ell+1}+\cdots +
(-1)^{\ell+n-2}i_n + (-1)^{\ell+n-1}
\frac{m}{2}+\lfloor\frac{\ell-1}{2}\rfloor
\end{eqnarray*}
and simplifying, we conclude
\begin{eqnarray*}
i_{\ell+1}\geq i_{\ell+2}-i_{\ell+3}+\cdots +
(-1)^{\ell+n} i_n + (-1)^{\ell+n+1}
\frac{m}{2}+\lfloor\frac{\ell+1}{2}\rfloor.
\end{eqnarray*}
\hfill{$\Box$}
\end{itemize}

   It is necessary to say a few words in order to evaluate the
novelty of  all this combinatorial
formulation. The following historical observations can be found in
[\cite{dieudonne89}; Ch. VI, Sect. 1.B, page 511]. In 1947, Steenrod
\cite{steenrod47} generalized the \v{C}ech-Whitney
definition  of the cup-product for a finite simplicial complex $K$
(in our context, we can consider $K$ as a polyhedral simplicial
set \cite{may67}).The idea was to keep several common vertices between
both factors
of the ``decomposition'' of a considered simplex, instead of one
as in the cup-product. In this way,
Steenrod established formulae for the cup-$i$ products, which were
``awkward to handle'', in his own words. He showed
that the
cohomology operations induced from this unwieldy description were
in fact independent of the order chosen on the vertices of $K$.

    Having this in mind, we could dare to say that in this paper
we are rediscovering  this old
description given by Steenrod and  clarifying it in a general
combinatorial framework. On the other hand, it is important to
note that we here determine the signs involved in the formulae of
the cup-$i$ products.

 Assuming that the face operators are evaluated in constant time,
the following result gives us a first measure of the computational complexity
of these formulae.

\begin{prp}
Let $\Z_2$ be the ground ring. Let $X$ be a simplicial set and $k$
a non-negative integer. If
$c\in C^{i+k}(X; \Z_2)$, then the number of face
operators taking part in the formula for $Sq^i(c)$ is $O(i^{k+1})$.
\end{prp}

\noindent{\bf Proof}

 Let $j=i+k$. Here, it is not necessary to distinguish the cases
 $i+j$ even and
 $i+j$ odd, since the proof is the same in both cases.

Firstly we count the number of summands
of the formula for $Sq^i(c)$ given in Corollary \ref{3}.

  The parameter $i_n$  contributes with
\begin{eqnarray}\label{rollo1}
m-S(n)+1=m-\lfloor\frac{m+1}{2}\rfloor - \lfloor
\frac{n}{2}\rfloor +1
\end{eqnarray}
summands in the formula, where $n=j-i$ and $m=i+j$.

Let us see that (\ref{rollo1}) is equal to $i+1$.
If $n$ is even then, this expression is equal to
$$m-\frac{m}{2}-\frac{n}{2}+1=\frac{m-n}{2}+1=i+1$$
and if $n$ is odd,
$$m-\frac{m+1}{2}-\frac{n-1}{2}+1=\frac{m-n}{2}+1=i+1.$$

  The parameter $i_{n-1}$ contributes with
$$i_n-S(n-1)=\lfloor\frac{m+1}{2}\rfloor - \lfloor
\frac{n-1}{2}\rfloor\leq i+1$$
summands.

In general,
the parameter $i_k$, with $1\leq k\leq
n-1$ contributes with
\begin{eqnarray}\label{rollo2}
i_{k+1}-S(k)=i_{k+2}-i_{k+3}+\cdots+(-1)^{k+n}
i_n+(-1)^{k+n+1}\lfloor\frac{m+1}{2}\rfloor - \lfloor
\frac{k}{2}\rfloor
\end{eqnarray}
summands, and using that $i_{\ell-1}-i_{\ell}\leq -1$ and $i_n\leq m$, it is
not difficult to see that the expression (\ref{rollo2}) is less or
equal to $i+1$ for each $1\leq k\leq n-1$.

Hence, the number
of summands appearing in the formulae of Corollary  \ref{3} is:
\begin{eqnarray*}
&&\sum_{i_n=S(n)}^m\,\cdots\,\sum_{i_3=S(3)}^{i_4-1}\;
\sum_{i_2=S(2)}^{i_3-1}i_2-S(1)\leq
\sum_{i_n=S(n)}^m\,\cdots\,\sum_{i_3=S(3)}^{i_4-1}\;
\sum_{i_2=S(2)}^{i_3-1}i+1\\\\
&&\leq(i+1)\sum_{i_n=S(n)}^m\,\cdots\,\sum_{i_3=S(3)}^{i_4-1}
i_3-S(2)\leq \cdots
\leq (i+1)^n.
\end{eqnarray*}

Since there are $m-n$ face operators in each summand and $n=k$,
the number of face operators that the formula has, is
less or equal to
$$(m-n)(i+1)^k= 2i(i+1)^k.$$
\hfill{$\Box$}

   We now discuss several facts about the
computation with the formulae of the Corollary \ref{3}. First of all,
it is clear that at least in the case in which  $X$ has
a finite number of non-degenerated simplices  in each
degree, our method can be seen as an actual algorithm for calculating
Steenrod squares.
For example, if the number of non-degenerated simplices in
every $X_{\ell}$ is $O({\ell }^2)$, then,
assuming that  each face operator of $X$ is an elementary operation,
the complexity of our
algorithm for calculating $Sq^i(c_{i+2})$ is $O(i^5)$. This
complexity
is obtained by the followings facts. On one hand, the number of
face operators taking part in the formula of
$Sq^{i}(c_{i+2})$ is $O(i^3)$. On the other hand, $Sq^i(c_{i+2})$ is a
$(2i+2)$-cochain and then, this cochain is
determined
by knowing its image over all the non-degenerated simplices in
$X_{2i+2}$.

Nevertheless, the most interesting examples appearing
in Algebraic Topology show, in general, a high
complexity in the number of non-degenerated simplices
in each degree. For example, let us take the
classifying space of a finite $2$-group. In this case,
the number of non-degenerated simplices in degree $\ell$
is $O(2^{\ell})$. Hence, our method will only be  useful
here in low dimensions. As we have mentioned in the introduction,
perhaps appropriately
combining these combinatorial formulae with classical
properties of Steenrod squares and well-known studies
for calculating cocycles  could allow us to
make a substantial  improvement in our method.

   Finally, this technique can be useful when dealing with
chain complexes arising from simplicial modules, like, for instance,
the Hochschild complex of an $R$-algebra $A$ is a simplicial
$R$-module (see, for example, \cite{Andre}). For
these differential graded modules, the use of
our simplicial method may be fruitful.

\section{A generalization to Steenrod reduced powers}

Real in \cite{real} established
formulae for the morphisms $\{D_i\}$ in terms of
the component morphisms of a given Eilenberg-Zilber
contraction. In an analogous way, we show here
that this result can be generalized to
Steenrod reduced powers (see \cite{steenrodepstein62,dieudonne89}).
It seems clear that this study must lead in the very near future
to an
explicit simplicial description of these cohomological operations.

Let us consider the Eilenberg-Zilber contraction:
$$(f,g,\phi): C_{*}^{\scst N}(X\times
\stackrel{\mbox{\scriptsize $p$ times}}\cdots \times X)
\ra C_{*}^{\scst N}(X)\ot \stackrel{\mbox{\scriptsize $p$
times}}\cdots \ot C_{*}^{\scst N}(X),$$
the diagonal $\Delta$ and the automorphisms $t$ and $T$ we have defined in
Section 2.

Now, the equality $tg=gT$ holds due to the associativity of $EML$
and to the good behavior of this morphism with regard to the
automorphisms $t$ and $T$.

Let us take a family of  automorphisms $\{\gamma_i \}_{i\geq 0}$
defined by
$$\gamma_{2j-1}=t\qquad \mbox{and}\qquad \gamma_{2j}=t+t^2+\cdots+t^{p-1}$$
and let
$\gamma=\gamma_i\gamma_{i-1}$.

Before proving the main result of this section, we need the following
propositions.

\begin{prp}\label{nose}
Let $p$ be  prime with $p\geq 2$ and $i$ a
non-negative integer, then
\begin{eqnarray*}
\phi d\gamma_i \cdots  \phi\gamma_{1}\phi &=&
(-1)^{i-1} \phi\gamma_i  \cdots
\phi\gamma_{1}d \phi+
\sum_{k=1}^{i-1} (-1)^{i-k}\phi \gamma_i  \cdots
\phi \gamma_{k+2}\phi \gamma\phi\gamma_{k-1}\cdots
\phi\gamma_{1}\phi.
\end{eqnarray*}
\end{prp}

\noindent {\bf Proof}

We prove the proposition by induction on the
parameter $i$.
\begin{itemize}
\item If $i=1$ then, $\phi d \gamma_1 \phi=\phi \gamma_1 d\phi$.
\item If $i=2$ then,
\begin{eqnarray*}
\phi d \gamma_2 \phi \gamma_1 \phi&=&
\phi \gamma_2 d\phi \gamma_1 \phi=
\phi \gamma_2 (gf-1-\phi d) \gamma_1 \phi\\
&=&-\phi \gamma \phi-\phi \gamma_2 \phi \gamma_1 d\phi.
\end{eqnarray*}
\item In general,
\begin{eqnarray*}
\phi d\gamma_i \phi \gamma_{i-1}\cdots  \phi\gamma_{1}\phi&=&
\phi\gamma_i d\phi \gamma_{i-1}\cdots  \phi\gamma_{1}\phi\\
&=&\phi \gamma_i(gf-1-\phi d)\gamma_{i-1}\cdots \phi\gamma_{1}\phi\\
&=&-\phi\gamma\phi\gamma_{i-2}\cdots  \phi\gamma_{1}\phi -
\phi\gamma_i \phi d\gamma_{i-1}\cdots  \phi\gamma_{1}\phi
\end{eqnarray*}
(by induction assumption)
\begin{eqnarray*}
&=&-\phi\gamma\phi\gamma_{i-2}\cdots  \phi\gamma_{1}\phi -
(-1)^{i-2} \phi \gamma_i \phi \gamma_{i-1}\cdots \phi \gamma_1
d\phi \\
&&-\sum_{k=1}^{i-2} (-1)^{i-1-k}
\phi\gamma_i \phi \gamma_{i-1}\cdots
\phi \gamma_{k+2}\phi \gamma\phi \gamma_{k-1} \cdots\phi\gamma_{1}\phi\\
&=&(-1)^{i-1} \phi \gamma_i \cdots \phi \gamma_1 d\phi\\
&&+\sum_{k=1}^{i-1} (-1)^{i-k}
\phi\gamma_i \cdots
\phi \gamma_{k+2}\phi \gamma\phi \gamma_{k-1}\cdots\phi\gamma_{1}\phi.
\end{eqnarray*}
\hfill{$\Box$}
\end{itemize}

Let  $\Gamma_i(k)=f\gamma_i \phi \gamma_{i-1}\cdots
\phi \gamma_{k+2}\phi \gamma\phi \gamma_{k-1}\phi
\cdots\phi\gamma_{1}\phi\Delta$. We can prove

\begin{prp}\label{util}
Let $p$ be prime and $i$ a non-negative integer, then
\begin{eqnarray*}
f\gamma_i \phi \gamma_{i-1}\cdots \phi \gamma_2\phi\Delta=
f\gamma_{i-1}\phi\gamma_{i-2} \cdots \phi\gamma_2\phi \gamma_1\phi \Delta
+\sum_{k=1}^{i-1}(-1)^{k+1}
\Gamma_i(k).
\end{eqnarray*}
\end{prp}

\noindent{\bf Proof}

Using $\gamma_{k+1}=\gamma_{k}+(-1)^{k+1} \gamma$,
for all $1<k<i$, we have
\begin{eqnarray*}
&&f\gamma_i \phi\gamma_{i-1}\cdots
\phi\gamma_{k+2}\phi \gamma_{k+1}\phi\gamma_{k-1}\cdots
\phi\gamma_1\phi\Delta\\
&&=f\gamma_i\phi \gamma_{i-1}\cdots\phi\gamma_{k+2}
\phi\gamma_k\phi\gamma_{k-1}\cdots \phi\gamma_1\phi\Delta+
(-1)^{k+1}\Gamma_i(k).
\end{eqnarray*}

If we use this fact successively, we obtain the following
identity:
\begin{eqnarray*}
f\gamma_i \phi \gamma_{i-1}\cdots \phi \gamma_2\phi \Delta =
f\gamma_{i-1}\phi\gamma_{i-2}\cdots\phi\gamma_2\phi\gamma_1\phi
\Delta+\sum_{k=1}^{i-1} (-1)^{k+1}\Gamma_i (k)
\end{eqnarray*}
\hfill{$\Box$}

The main result of this section is the following one.

\begin{thr}
Let $p$ be a prime number, $p\geq 2$, and $i$ a
non-negative integer. Then, there exists a sequence of
morphisms $\{D_i\}$ defined by
$$D_i=f \gamma_i \phi
\gamma_{i-1}\cdots  \phi\gamma_{1} \phi \Delta,$$
verifying that
$$d D_i + (-1)^{i+1} D_i d=\alpha_i D_{i-1},$$
where $\alpha_{2j-1}=T-1$ and
$\alpha_{2j}=1+T+T^2+\cdots+T^{p-1}$.
\end{thr}

\noindent{\bf Proof}

Let us begin with the first term of the identity:
\begin{eqnarray*}
d D_i+(-1)^{i+1}D_i d &=&
d f \gamma_i \phi\gamma_{i-1}\cdots  \phi\gamma_{1}\phi \Delta +
(-1)^{i+1}f\gamma_i\phi\gamma_{i-1}\cdots \phi\gamma_{1}\phi \Delta
d\\
&=& f \gamma_i d \phi\gamma_{i-1}\cdots  \phi\gamma_{1}\phi \Delta \\
&&+(-1)^{i+1}f\gamma_i\phi\gamma_{i-1}\cdots \phi\gamma_{1}
(gf-1-d\phi)\Delta\\
&=&f \gamma_i d\phi \gamma_{i-1}\cdots  \phi\gamma_{1}\phi \Delta\\
&&+(-1)^i f \gamma_i \phi\gamma_{i-1}\cdots  \phi\gamma_{2}\phi
\Delta+
(-1)^i f\gamma_i \phi \gamma_{i-1}\cdots \phi
\gamma_{1}d\phi \Delta
\end{eqnarray*}
(by Proposition \ref{nose}, we get)
\begin{eqnarray*}
&=& f\gamma_i
d\phi\gamma_{i-1}\cdots \phi
\gamma_1\phi\Delta+(-1)^i f\gamma_i\phi\gamma_{i-1} \cdots \phi
\gamma_2\phi \Delta\\&&+(-1)^i(-1)^{i-2}f\gamma_i \phi
d \gamma_{i-1}\cdots \phi \gamma_1 \phi
\Delta\\
&&+\sum_{k=1}^{i-2} (-1)^i(-1)^kf\gamma_i \phi
\gamma_{i-1}\cdots  \phi \gamma_{k+2}\phi
\gamma\phi \gamma_{k-1}\phi
\cdots\phi\gamma_{1}\phi\Delta\\
&=&f\gamma_i (d\phi+\phi
d)\gamma_{i-1}\cdots \phi \gamma_1\phi\Delta\\
&&+(-1)^i
f\gamma_i \phi\gamma_{i-1}\cdots \phi \gamma_2\phi
\Delta+\sum_{k=1}^{i-2} (-1)^{i-k} \Gamma_i(k)\\
&=& f\gamma_i
(gf-1)\gamma_{i-1}\cdots \phi
\gamma_1\phi\Delta\\
&&+(-1)^i f\gamma_i \phi\gamma_{i-1}\cdots \phi
\gamma_2\phi \Delta+\sum_{k=1}^{i-2}
(-1)^{i-k}\Gamma_i(k)
\end{eqnarray*}
(let
$\beta_{2j}=T+T^2+\cdots+T^{p-1}$ and $\beta_{2j-1}=T$,
then we have)
\begin{eqnarray*}
&=& \beta_i f
\gamma_{i-1}\phi\gamma_{i-2}\cdots \phi \gamma_1\phi\Delta\\
&&+(-1)^i
f\gamma_i\phi\gamma_{i-1} \cdots \phi \gamma_2\phi
\Delta+\sum_{k=1}^{i-1}
(-1)^{i-k} \Gamma_i(k)
\end{eqnarray*}
(Proposition \ref{util} implies)
\begin{eqnarray*}
&=&\beta_i f
\gamma_{i-1}\phi\gamma_{i-2}\cdots \phi \gamma_1\phi\Delta+(-1)^i
f\gamma_{i-1}\phi\gamma_{i-2} \cdots \phi \gamma_1\phi
\Delta\\
&&+\sum_{k=1}^{i-1}
(-1)^{i+k+1} \Gamma_i(k)+
\sum_{k=1}^{i-1}
(-1)^{i-k} \Gamma_i(k)\\
&=&\alpha_i f
\gamma_{i-1}\phi\gamma_{i-2}\cdots \phi \gamma_1\phi\Delta.
\end{eqnarray*}
\hfill{$\Box$}

\section{Proof of the main theorem }

The proof consists in finding out the factors  of the
formula (written
in the canonical way) that are degenerated
and in eliminating the summands having
these factors.

First of all, notice that
using  the commutativity properties of the operators of a simplicial set
(essentially, {\bf (s3)}),
it is easy to see that a factor of the formula
whose expression begins (on the left) by
\begin{eqnarray}\label{formadeg}
\partial_{j_1} \cdots \partial_{j_t}s_{k}\cdots
\end{eqnarray}
such that $0\leq j_1<\cdots <j_t<k$,
is degenerated in its simplified form.

Having said that, let us begin with the proof of the
theorem.

For $n=0$, we obtain the explicit formula for the
Alexander-Whitney operator.

Let us assume that the formula is
true for $k\leq n$ so, let us prove that the formula is
true for the case $n+1$.

Let us consider that $n$ is even (when $n$ is odd, the
proof is similar). In this case, by induction
assumption, the formula over an element of degree $m$ is as follows:

$$\begin{array}{ll}
AW\;(t\; SHI)^{n+1}\;=&AW\;(t\; SHI)^n\;(t\; SHI)\\\\
=\displaystyle\sum_{i_n=n}^{m+1}\,\cdots
\,\displaystyle\sum_{i_0=0}^{i_1-1}\quad
\displaystyle\sum&
(-1)^{A(n)+B(n,m+1,\bar{\i})+C(n,m+1,\bar{\i})+D(n,\bar{\i})}
(-1)^{\bar{m}+sig(\alpha,\beta)+1}\\
&\partial_{i_0+1}\cdots\cdot
\partial_{m+1}s_{\alpha _{p+1}+\bar{m}} \cdots
s_{\alpha _1+\bar{m}} \partial _{\bar{m}} \cdots
\partial _{m-q-1} \\&\ot \partial_0 \cdots\cdot
\partial_{i_{n}-1}s_{\beta _q+\bar{m}} \cdots s_{\beta_1+\bar{m}}
s_{\bar{m}-1} \partial _{m-q+1} \cdots
\partial _m
\end{array}$$

\noindent where $\bar{\i}=(i_0,i_1,\dots,i_n)$, $\bar{m}=m-p-q$,
$sig(\alpha,\beta)=\displaystyle\sum_{i=1}^{p+1}
\alpha_i-(i-1)$, and the last sum is
taken over all the indices $0 \leq q \leq m-1$, $0
\leq p \leq m-q-1$ and $( \alpha , \beta ) \in \{ (p+1,q)$-shuffles$\}$.

Let us recall that if  the formula is in the
normalized form, the summands  which have a
degenerated factor  must be
eliminated. The following cases are considered:

\noindent If $i_n>m-p$, we have to consider the following cases:

\begin{itemize}
\item If $i_n=m-p+t$ and
$\beta_q< q-1+t$, with
$1\leq t\leq p$ then,
$$\alpha_{p+1}=p+q>\cdots
>\alpha_{t+1}=q+t>\alpha_t=q+t-1>\beta_q.$$
So, the first factor of these summands is:
$$\partial_{i_0+1}\cdots \cdot
\partial_{i_{n-1}-1}\partial_{m-p+t+1}\cdots
\partial_{m+1}s_m\cdots
s_{m-p+t-1}s_{\alpha_{t-1}}\cdots
s_{\alpha_1}\partial_{\bar{m}}\cdots \partial_{m-q-1}$$
$$=\partial_{i_0+1}\cdots \cdot
\partial_{i_{n-1}-1}s_{m-p+t-1}s_{\alpha_{t-1}}\cdots
s_{\alpha_1}\partial_{\bar{m}}\cdots \partial_{m-q-1}.$$
Since $i_{n-1}-1<i_n-1=m-p+t-1$, this factor is degenerated by
(\ref{formadeg}).

\item If $i_n=m-p+t$ and $\beta_q=q-1+t$, with
$1\leq t\leq p$ then,
$$\alpha_{p+1}=p+q>\cdots >\alpha_{t+1}=q+t>\beta_q=q+t-1.$$
And hence, in this case, the expression of the first
factor begins in the form:
$$\partial_{i_0+1}\cdots
\cdot\partial_{i_{n-1}-1}\partial_{m-p+t+1}\cdots
\partial_{m+1}s_m\cdots
s_{m-p+t}s_{\alpha_t+\bar{m}}\cdots$$
$$=\partial_{i_0+1}\cdots
\cdot\partial_{i_{n-1}-1}s_{\alpha_t+\bar{m}}\cdots.$$

\newpage

Now, we have to consider two different cases:
\begin{itemize}
\item If $i_{n-1}-1<\alpha_t+\bar{m}$ then, this factor
is degenerated.

\item If $i_{n-1}-1\geq \alpha_t+\bar{m}$, let us
denote $\alpha_t=q+t-1-j$, where $1\leq j\leq q+t-1$.
Then, $$\alpha_{p+1}=p+q>\cdots >
\alpha_{t+1}=q+t,$$
and
$$\beta_q=q+t-1>\cdots >
\beta_{q-j+1}=q+t-j>\alpha_t=q+t-1-j. $$
Hence, the second factor of these
summands is:
$$\partial_0\cdots \cdot \partial_{i_{n-2}-1}
\partial_{i_{n-1}+1}\cdots \partial_{m-p+t-1}
s_{m-p+t-1}\cdots
s_{m-p+t-j}s_{\beta_{q-j}+\bar{m}}\cdots$$
$$=\partial_0\cdots \cdot \partial_{i_{n-2}-1}
s_{i_{n-1}}\cdots s_{m-p+t-j}s_{\beta_{q-j}+\bar{m}}\cdots
$$
Since
$i_{n-2}-1<i_{n-1}$, this factor is degenerated.
\end{itemize}

\item If $i_n=m-p+t$ and $\beta_q>q-1+t$, with $1\leq t\leq p$ then
$i_n-1<\beta_q+\bar{m}.$
So, the second factor of the summands has
the form (\ref{formadeg}) and hence, these summands must be
eliminated.

\item If $i_n=m+1$ and $\beta_q<p+q$ then $\alpha_{p+1}=p+q$ and
since $i_{n-1}-1<i_n-1=m$ then, the first factor of the
summands is degenerated.

\item If  $i_n=m+1$ and $\beta_q=p+q$ then
$$\beta_q=p+q>\cdots >\beta_{j+1}=p+j+1>
\alpha_{p+1}=p+j$$
with $0\leq j\leq q-1$ and the first
factor of the summands is:
$$\partial_{i_0+1}\cdots\cdot
\partial_{i_{n-1}-1}s_{m-q+j}
s_{\alpha_{p}+\bar{m}}\cdots
s_{\alpha_1+\bar{m}} \partial_{\bar{m}} \cdots
\partial_{m-q-1}.$$

We have to consider two different cases:
\begin{itemize}
\item If $i_{n-1}-1<m-q+j$ then, this factor is
degenerated.
\item If $i_{n-1}-1\geq m-q+j$ then
the second factor of the summands is:
$$\partial_0\cdots \cdot
\partial_{i_{n-2}-1}\partial_{i_{n-1}+1}\cdots
\partial_m  s_m\cdots
s_{m-q+j+1}s_{\beta_j+\bar{m}}\cdots$$
$$=\partial_0\cdots \cdot \partial_{i_{n-2}-1}
s_{i_{n-1}}\cdots s_{m-q+j+1}
s_{\beta_j+\bar{m}}\cdots,$$
which is degenerated
due to the fact that $i_{n-2}-1<i_{n-1}.$
\end{itemize}
\end{itemize}

\noindent If $i_n<m-p$, then
$i_n-1<\beta_q+\bar{m}.$
So, these summands have the second
factor in the form (\ref{formadeg}) and hence, must be
eliminated.

\noindent If $i_n=m-p$, two cases hold:
\begin{itemize}
\item If $\beta_q>q-1$, then the second
factor of these  summands is degenerated as above.

\item If $\beta_q=q-1$
and $i_{n-1}>\bar{m}-2$ then
$$\alpha_{p+1}=p+q>\cdots>\alpha_1=q>\beta_q=q-1>\cdots>\beta_1=0,$$
and the second factor of the tensor product is
$$\partial_{0}\cdots \cdot \partial_{i_{n-2}-1}
\partial_{i_{n-1}+1} \cdots \partial_{m-p-1}s_{m-p-1} \cdots
s_{m-p-q-1}\partial _{m-q+1}
\cdots \partial _m$$
$$=\partial_{0}\cdots \cdot \partial_{i_{n-2}-1}
s_{i_{n-1}}\cdots s_{m-p-q-1}\partial _{m-q+1}
\cdots \partial _m;$$

\noindent and since $i_{n-2}-1<i_{n-1}$,
this factor is degenerated.

Finally, if $\beta_q=q-1$ and $i_{n-1}\leq
\bar{m}-2$
then the formula  (save for the signs) corresponding to
$AW\;(t\;SHI)^{n+1}$ is:

\begin{eqnarray*}
&&\sum_{i_{n-1}=n-1}^{\bar{m}-2}\,
\cdots \, \sum_{i_0=0}^{i_1-1}\quad
\sum_{q=0}^{m-1}\quad \sum_{p=0}^{m-q-1}\\\\
&&\quad \partial_{i_0+1}\cdots \cdot
\partial_{i_{n-1}-1} \partial_{m-p+1}
\cdots \partial_{m+1}s_m \cdots s_{m-p} \partial _{\bar{m}}
\cdots \partial_{m-q-1} \\
&& \quad \ot \partial_0 \cdots \cdot
\partial_{i_{n-2}-1}
\partial_{i_{n-1}+1}
\cdots \partial_{m-p-1}s_{m-p-1} \cdots s_{\bar{m}-1}
\partial _{m-q+1} \cdots \partial _m   \\\\
&&= \sum_{i'_{n+1}=n+1}^{m}\quad
\sum_{i'_n=n}^{i'_{n+1}-1}\quad
\sum_{i_{n-1}=n-1}^{i'_{n}-1}\,\cdots \,
\sum_{i_0=0}^{i_1-1}\\\\
&& \quad \partial_{i_0+1}\cdots \cdot
\partial_{i_{n-1}-1}  \partial _{i'_n+1} \cdots \partial
_{i'_{n+1}-1}
 \\&& \quad \ot \partial_0 \cdots \cdot
\partial_{i_{n-2}-1}\partial_{i_{n-1}+1}\cdots
\partial _{i'_n-1}\partial _{i'_{n+1}+1}
\cdots \partial_{m};\end{eqnarray*}

\noindent where $i'_n=\bar{m}-1$ and $i'_{n+1}=m-q$.
\end{itemize}

Now, let us study the signs of the formulae in this last case.
Keeping in mind that we are working with the exponent of
$(-1)$, all the identities are mod $2$.

This proof is based on
the fact that if and only if $i_n=m-p$ and $\beta_q=q-1$,
the summands of $AW\;(t\;SHI)^n$ are non-degenerated.

We first verify the formula in the case $n=1$.
 The exponent of $(-1)$ associated to each
summand is:
\begin{eqnarray*}
&&\bar{m}+sig(\alpha,\beta)+1=
\bar{m}+\sum_{i=1}^{p+1}\alpha_i-(i-1)+1\\
&&=\bar{m}+\sum_{i=1}^{p+1}(q-1+i-(i-1))+1= \bar{m}+q(p+1)+1
\end{eqnarray*}
(note $i_0=\bar{m}-1$ and $i_1=m-q$)
$$=i_0+(m+i_1)(i_0+i_1)=A(1)+B(1,m,\bar{\i})+C(1,\bar{\i})+D(1,m,\bar{\i})$$
where $\bar{\i}=(i_0,i_1)$, as asserted.

In general, we have to prove that
$$A(n)+B(n,m+1,\bar{\i})+C(n,\bar{\i})+D(n,m+1,\bar{\i})+
\bar{m}+sig(\alpha,\beta)+1$$
$$=A(n+1)+B(n+1,m,\bar{\i})+C(n+1,\bar{\i})+D(n+1,m,\bar{\i}).$$

\begin{itemize}
\item If $n$ is even then
$D(n,m+1,\bar{\i})=0$ and the exponent of $(-1)$ in each summand is:
\begin{eqnarray*}
&&A(n)+B(n,m+1,\bar{\i})+C(n,\bar{\i})+\bar{m}+sig(\alpha,\beta)+1\\
&&=A(n)+B(n,m+1,\bar{\i})+C(n-1,\bar{\i})\\
&&\quad+(m+p+i_{n-1})(i_{n-1}+\cdots +i_0)+
\bar{m}+q(p+1)+1
\end{eqnarray*}
(since $i'_n=\bar{m}-1$, $i'_{n+1}= m-q$ and these identities are mod
$2$ then)
\begin{eqnarray*}
&&=A(n)+B(n,m+1,\bar{\i})+C(n-1,\bar{\i})\\
&&\quad+(m+1+i'_{n+1}+i'_n+i_{n-1})(i_{n-1}+\cdots +i_0) \\
&&\quad+i'_n+(m+i'_{n+1})(i'_{n+1}+i'_n)\\
&&=A(n)+B(n,m+1,\bar{\i})+C(n-1,\bar{\i})+i_{n-1}+\cdots+i_0\\
&&\quad+(i'_n+i_{n-1})(i_{n-1}+\cdots+i_0)
+(m+i'_{n+1})(i_{n-1}+\cdots+i_0)\\
&&\quad+i'_n+(m+i'_{n+1})(i'_{n+1}+i'_n)\\
&&=A(n)+B(n,m+1,\bar{\i})+i'_n+i_{n-1}+\cdots+i_0\\
&&\quad+C(n+1,\bar{\i})+D(n+1,m,\bar{\i}).
\end{eqnarray*}

\newpage

We have to distinguish two cases:
\begin{itemize}
\item If  $n\equiv 0$ mod $4$ then
$A(n)=A(n+1)$ and
\begin{eqnarray*}
&&A(n)+B(n,m+1,\bar{\i})+i'_n+i_{n-1}+\cdots+i_0\\\\
&&=A(n+1)+\sum_{j=0}^{\frac{n-2}{2}}i_{2j+1}+i'_n+i_{n-1}+\cdots+i_0\\\\
&&=A(n+1)+\sum_{j=0}^{\frac{n-2}{2}}i_{2j}+i'_n=A(n+1)+B(n+1,m,\bar{\i}).
\end{eqnarray*}

\item If $n\equiv 2$ mod $4$ then
$A(n)=A(n+1)+1$ and
\begin{eqnarray*}
&&A(n)+B(n,m+1,\bar{\i})+i'_n+i_{n-1}+\cdots+i_0\\\\
&&=A(n+1)+1+\sum_{j=0}^{\frac{n-2}{2}}i_{2j}+
m+p+i'_n+i_{n-1}+\cdots+i_0\\\\
&&=A(n+1)+1+\sum_{j=0}^{\frac{n-2}{2}}i_{2j}+
m+i'_{n+1}+i'_n+1\\\\
&&\quad+i'_n+i_{n-1}+\cdots+i_0\\\\
&&=A(n+1)+\sum_{j=0}^{\frac{n-2}{2}}i_{2j+1}+i'_{n+1}+m=
A(n+1)+B(n+1,m,\bar{\i}).
\end{eqnarray*}
\end{itemize}

\item If $n$ is odd
then the exponent of $(-1)$ is:
\begin{eqnarray*}
&&A(n)+B(n,m+1,\bar{\i})+C(n,\bar{\i})+D(n,m+1,\bar{\i})+
\bar{m}+sig(\alpha,\beta)+1\\
&&=A(n)+B(n,m+1,\bar{\i})+C(n,\bar{\i})\\
&&\quad+(m+1+m+p)(m+p+i_{n-1}+\cdots+i_0)+\bar{m}+q(p+1)+1
\end{eqnarray*}
  (since $i'_n=\bar{m}-1$, $i'_{n+1}= m-q$ and these identities are mod
$2$ then)
\begin{eqnarray*}
&&=A(n)+B(n,m+1,\bar{\i})+C(n,\bar{\i})\\
&&\quad+(i'_{n+1}+i'_n)(m+1+i'_{n+1}+i'_n+i_{n-1}+\cdots +i_0)\\
&&\quad+i'_n+(m+i'_{n+1})(i'_{n+1}+i'_n)\\
&&=A(n)+B(n,m+1,\bar{\i})+C(n,\bar{\i})+
(i'_{n+1}+i'_n)(i'_n+i_{n-1}+\cdots+i_0)\\
&&\quad+i'_{n+1}+i'_n+(i'_{n+1}+i'_n)(m+i'_{n+1})+
i'_n+(m+i'_{n+1})(i'_{n+1}+i'_n)\\
&&=A(n)+B(n,m+1,\bar{\i})+C(n,\bar{\i})\\
&&\quad+(i'_{n+1}+i'_n)(i'_n+i_{n-1}+\cdots+i_0)+i'_{n+1}\\
&&=A(n)+B(n,m+1,\bar{\i})+i'_{n+1}+C(n+1,\bar{\i})+D(n+1,m,\bar{\i}).
\end{eqnarray*}

Since $n$ is odd then $n\equiv 1,3,5,7$ mod $8$  and, in these cases,
$A(n)=A(n+1)$.

Now, we have to distinguish two cases:
\begin{itemize}
\item If  $n\equiv 1$ mod $4$ then
\begin{eqnarray*}
B(n,m+1,\bar{\i})+i'_{n+1} =\sum_{j=0}^{\frac{n-1}{2}}i_{2j}+i'_{n+1}=
\sum_{j=0}^{\frac{n+1}{2}}i_{2j}=B(n+1,m,\bar{\i}).
\end{eqnarray*}

\item If $n\equiv 3$ mod $4$ then
\begin{eqnarray*}
B(n,m+1,\bar{\i})+i'_{n+1}&=&
\sum_{j=0}^{\frac{n-3}{2}}i_{2j+1}+i_n+m+1+i'_{n+1}\\\\
&=&\sum_{j=0}^{\frac{n-3}{2}}i_{2j+1}+i'_n+m+1+m+1\\\\
&=&\sum_{j=0}^{\frac{n-1}{2}}i_{2j+1}=B(n+1,m,\bar{\i}).
\end{eqnarray*}
\end{itemize}
\hfill{$\Box$}
\end{itemize}

\end{document}